\documentclass[graybox]{svmult}

\usepackage{type1cm}        % activate if the above 3 fonts are
\usepackage{makeidx}         % allows index generation
\usepackage{graphicx}        % standard LaTeX graphics tool
\usepackage{multicol}        % used for the two-column index
\usepackage[bottom]{footmisc}% places footnotes at page bottom

\usepackage{newtxtext}       %
\usepackage{newtxmath}       % selects Times Roman as basic font

\makeindex             % used for the subject index
\newcommand{\be}{\begin{equation}}
\newcommand{\ee}{\end{equation}}
\newcommand{\R}{\mathbb{R}}
\newcommand{\vphi }{\varphi }
\begin{document}

\title*{On a class of self-similar solutions of the Boltzmann equation}
\author{ A.V. Bobylev}
\institute{ A.V. Bobylev \at Keldysh Institute of Applied Mathematics RAS, 125047 Moscow, RUSSIA, \email{alexander.bobylev@kau.se}
}
\maketitle
\begin{abstract}
{We consider a class of distribution functions
 having the form $ f(v,t)=t^{-dat} F(v e^{-at} )$,
 $ a=\mathrm{const.} $, where $ v \in \R^{^{d}}, d \geq 2 $
 and $ t \in \R_{+} $ denote the particle velocity and
 the time. This class of self-similar solutions
 to the spatially homogeneous Boltzmann equation (BE) for
 Maxwell molecules was studied by Bobylev and Cercignani in
 early 2000s. The solutions are positive, but have an
 infinite second moment (energy).
 % Therefore they are not
 %very interesting for applications.
  However, the same class of distribution functions with finite energy appears
 to be closely connected with quite different class of
 group-invariant solutions of the spatially inhomogeneous
 BE. This is a motivation for considering the so-called
 modified spatially homogeneous BE, which contains an extra
 force term proportional to a given matrix $ A $. The  modified
 BE was recently studied under assumption of "sufficient
 smallness of norm $ \|A\| $" without explicit estimates of
 the smallness. We fill this gap and prove that all
 important facts related to self-similar solutions remain
 valid also for $ \|A\| =O[10^{-1}]$ in appropriate dimensionless units.}
\end{abstract}
\section{Introduction}
\label{sec:1}
\numberwithin{equation} {section}
Self-similar solutions of nonlinear equations of mathematical physics are always interesting for both physicists and mathematicians. The simplest example of such solution for the classical Boltzmann equation can be described as follows. We consider the spatially homogeneous case of this equation for Maxwell molecules (with or without cut-off). The solution
$f(v, t)$ of this equation depends on the velocity
$ v \in \R^3 $ and the time $t \in \R_{+} $.
It is easy to show that the Boltzmann equation admits (at least formally) a particular class of solutions having the form
\begin{equation}\label{1.1}
  f(v, t) = e^{-3ct} F(v e^{-ct}), \quad c = const.                                       \end{equation}
Then we obtain the time-independent equation for the self-similar profile $F(v)$ and some conditions on the constant $c$. This is what we did jointly with Carlo Cercignani approximately 20 years ago~\cite{BC1, BC2, BC3}. By that time I already had some experience with related classes of group – invariant solutions to this equation~\cite{B3,B4}. On the other hand, Carlo Cercignani had an idea to use the Pomea~\cite{Pomea} approach to the structure of infinitely strong shock wave. The connection with this approach is briefly explained in~\cite{BC2}. Unfortunately we did not manage to produce something useful for shock waves. However, the results of~\cite{BC1, BC2, BC3} for the spatially homogeneous case were interesting. In fact we studied only isotropic (depending on $|v |$) solutions having the form
\eqref{1.1}. We have proved that these solutions are positive and describe a large–time asymptotics for certain classes of initial data. Integral representations of two particular solutions were constructed in explicit form. They represent two non-trivial examples of eternal solutions of the Boltzmann equations, which exist for all real values of  time $t$ (for more information see the book~\cite{Bobook}).

The self-similar solutions of the form~\eqref{1.1} and related questions were studied more deeply in some interesting mathematical papers, see,
 e.g.~\cite{Can.}, \cite{Mor.}. Of course, a serious drawback of solutions~\eqref{1.1}  is that they obviously contradict to the energy (the second moment of $F$) conservation. The contradiction can be avoided only if  the second moment of the positive function $F(v)$ is infinite. It is really so for above discussed solutions, they all have an infinite second moment. Thus, it is a very unusual for rarefied gas dynamics class of solutions.
However, there is a quite different class of spatially inhomogeneous solutions to the Boltzmann equation, which can also lead to distribution functions of the form (1.1). These are the so-called homoenergetic affine flows introduced in 1950s independently by Galkin~\cite{GA} and Truesdell~\cite{Tru.}.
There are many related references, we just mention two
books~\cite{Santos}, \cite{TruMun},
 and contributions of Cercignani to these areas~\cite{Cer1},
 \cite{Cer2}, \cite{Cer3}. Recently these flows were considered in detail in series of papers by James, Nota and Velazquez~\cite{JNV}, \cite{JNV1}, \cite{JNV2}, see also~\cite{BNV}, \cite{Duan}. Roughly speaking, these  are spatially inhomogeneous flows of  gas having the linear (with respect to spatial variable $x \in \mathbb{R}^3$ ) profile of the bulk velocity. This assumption allows to reduce the problem to the modified spatially homogeneous Boltzmann equation. The general form of this equation introduced in~\cite{BNV} is also considered in the present paper. Our aim is to study solutions of this equation having the form~\eqref{1.1}. The interest to this problem is partly caused by the first proof  of existence of self-similar solutions of that kind in~\cite{JNV}  and~\cite{BNV}.  The modified Boltzmann equation depends on a small parameter. This small parameter is assumed  in these papers to be "as small as we want".  One of the goals of the present paper is to weaken this assumption and to prove that all main results of our previous work~\cite{BNV} remain valid for moderately small values of the parameter (roughly $10\%  $  of contribution of collisions) and to show why the proof cannot be extended to larger values. Unfortunately the restricted volume of the paper does not allow to include details of the proofs, they will be published elsewhere.  Below we confine ourselves to statement of the problem in Sect.~\ref{sec:2},  eigenvalue problem for matrices in
Sect.~\ref{sec:3} and formulation of  main results in Sect.~\ref{sec:4}.
\section{Statement of the problem}
\label{sec:2}
We consider  a modified Boltzmann equation~\cite{BNV} for the
distribution function $f(v,t)$, where  $v\in \R^d$,
 $d\geq 2$, and $t\in \R_{+}$ denotes respectively the particle
 velocity and time. The equation reads
\be \label{2.1} %   1.1
\partial_{t}f- div_{v}\left(A v f\right)    =Q\left(f,f\right),
\ee
where $ A\in M_{d\times d}(\R) $ is a constant  real matrix and
 $ Q\left( f,f \right) $ is the collision integral for Maxwell
 molecules
\be \label{2.2} %1.2
Q\left(f,f\right) \left(  v\right)= \int \limits_{\mathbb{R}^{d}
\times S^{d-1}} \!\! \!dw dn  g\left(  \hat{u} \cdot n \right)\left[  f(v^{\prime})f(w^{\prime})- f(v) f(w)\right],
\; n\in S^{d-1}
\ee
\[
u=v-w,\quad
\hat{u}={u} /   {\vert u \vert },\quad
v^{\prime} =\frac{1}{2} \left(v+w+\vert u \vert n \right), \;
w^{\prime}=\frac{1}{2}\left(v+w-\vert u \vert n \right).
\]
We assume that
\be \label{2.3} %1.3
f(v,0)=f_0(v)\geq 0, \quad \int_{\R^d} dv \ f_0(v) v=0,
 \quad \int_{\R^d} dv \ f_0(v)=1.
\ee
The kernel $ g(\eta) $ in~\eqref{2.2}, with $ \eta \in [-1,1] $,
is non-negative and normalized by unity
\be\label{2.4} %1.4
\int_{S^{d-1}} dn g(\omega \cdot n)=1,
\quad
\omega \in S^{d-1}.
\ee
The motivation for considering the modified Boltzmann equation~\eqref{2.1}
is discussed in detail in~\cite{BNV}. In particular, we note
that the well-known  shear flow for the Boltzmann equation is
described by equation~\eqref{2.1} with nilpotent matrix $ A $
such that
\begin{equation}\label{2.5} %2.1
A=\{a_{ij}; \quad i,j=1,...,d\},
\end{equation}
where $ a_{12}=a=\mathrm{const.} $ and all other elements of
$ A $ are zeros. Many related details and references can be
found in~\cite{TruMun}, \cite{Santos},
\cite{Cer1}, \cite{Cer2}, \cite{Cer3},
  \cite{Bob.- Spi.}, \cite{JNV}.

Formally the second term in~\eqref{2.1} describes the action
of the external force
\be \label{2.6}
F=-Av, \quad  A\in M_{d\times d}(\R),
\ee
which looks like the anisotropic \emph{friction} force proportional
to components of the particle velocity. Let us consider e.g., the simplest case of equation~\eqref{2.1} with
\be \label{2.7}
A=a I, \quad  a\in \R,
\ee
where $ I $ is the unit matrix, $ a $ is a constant with any
sign. If $ a > 0 $ this is just a regular friction force
$ F=- a v $, $ a > 0 $. The solution $ f(v,t) $
 of~\eqref{2.1}--\eqref{2.4} under assumption~\eqref{2.7}
 leads to the following behaviour of the second moment (energy):
 \begin{equation*}
 {d\mathcal{E}(t)}/{dt}= - 2 a \mathcal{E}(t)
 \;
 \Rightarrow\;
 \mathcal{E}(t)=
 \mathcal{E}(0)e^{-2 a t},
 \end{equation*}{where}
\[\quad
\mathcal{E}(t)=\int_{\R^d} dv f(v,t)|v|^2.
\]
This equality gives an idea to consider the equation~\eqref{2.1}
with $ Av $ in self-similar variables by substitution
\be\label{2.8}
f(v,t)=e^{d a t} \tilde{f}(\tilde{v},t),
\quad
\tilde{v} = v e^{a t}.
\ee
 Then after simple calculations, we obtain the familiar spatially
homogeneous Boltzamnn equation for $ \tilde{f}(\tilde{v},t) $
\be
\label{2.9}
\hat{f}_{t}= Q(\tilde{f},\tilde{f}), \quad
\tilde{f}|_{t=0} =f_0(\tilde{v}).
\ee
If, in addition,
\begin{equation*}
 \mathcal{E}(0)= \int_{\R^d} dv f_{0}(v)|v|^2=d,
\end{equation*}
 then we know (H-theorem for~\eqref{2.9}) that
 \begin{equation*}
 \tilde{f}(\tilde{v},t)
 \;  \rightarrow\;
 \tilde{f}_{M}(\tilde{v})=
 (2\pi)^{-d/2 } e^{- |\tilde{v}|^2 /2 },
 \quad
 \text{as} \;t \to 0.
 \end{equation*}
 Hence, coming back to initial variables, we
(1)
 obtain the
 simplest self-similar solution of~\eqref{2.1}, \eqref{2.7},
 namely
\begin{equation}
\label{2.10}
 f_{s-s}(v,t)=\left(2\pi e^{-2 a t} \right)^{-d/2}
 \exp \left( - {|v e^{ a t}|^{2}}/{2}
 \right),
  \end{equation}
and
(2) show that this particular solution is an attractor for
various classes of initial data.
It is obvious, that this simple example is valid also for
$ a <0 $  in~\eqref{2.7} (accelerating forces)  and for arbitrary
kernel (not necessary the Maxwellian one) in the collision
integral.

Roughly speaking, our task is to prove that the situation is,
to some extent, similar in the case of arbitrary matrix
$ A $ in~\eqref{2.1} provided that its norm is not too large.
In fact, all proofs were already done in our previous
paper~\cite{BNV} with standard formulations of results like as
"There exists such $ \varepsilon_{0} > 0 $ that the following
property holds under assumption that
 $ \|A\| \le \varepsilon_{0}  $... ". This approach allows to
 avoid some technical work, but it does not show true limits
 (in terms of $ \|A\| $) of the results. The main aim of this paper
 is to partly clarify this question. Here and below we use
 the so-called operator norm for matrices~\cite{Kato}.
 Its properties are discussed in the next section.

 Following~\cite{BNV} we pass to  the Fourier-representation~\cite{B75} of
 the equation~\eqref{2.1} and introduce the characteristic function~\cite{Feller}
 $\varphi(k,t)$
\be \label{2.11}
\varphi(k,t)=\int_{\R^d} dv \ f(v,t)e^{-i k\cdot v}, \;\; k\in\R^d.
\ee
Then we obtain
\be \label{2.12}
\partial_t \varphi  + \left(Ak\right)\cdot \partial_k \varphi =
 \mathcal{I}^{+}(\varphi,\varphi)-\varphi_{|_{k=0}} \varphi\,,
\ee
where
\begin{align}  \label{2.13}
%\begin{split}
 \mathcal{I}^{+}(\varphi,\varphi)(k)= \!\int \limits_{S^{d-1}}
 \!\!dn
  g\left(  \hat{k} \cdot n \right)\varphi(k_{+}) \varphi (k_{-}),
 \;
 % \\
\; k_{\pm}=\frac{1}{2}\left( k\pm \vert k \vert n\right),
 \;
 \hat{k}=\frac{k}{\vert k\vert}.
 %\end{split}
\end{align}

The initial condition becomes
\begin{equation*}
\varphi (k,0)=\varphi_0(k)= \int_{\R^d} dv
 \ f_0(v)e^{-i k\cdot v}, \quad\varphi_0(0)=1.
\end{equation*}
Note that  \eqref{2.1} implies the mass conservation. Therefore
\be \label{2.14}
\varphi(0,t)=\varphi_0(0)=1,
\ee
and we obtain from~\eqref{2.12}
\be \label{2.15}
\partial_t \varphi  +\varphi + \left(Ak\right)\cdot \partial_k \varphi
 = \mathcal{I}^{+}(\varphi,\varphi) = \Gamma(\varphi).
\ee
For brevity we consider below the self-similar solution only.
Following~\cite{BNV}, we look for such solution in the form
\be \label{2.16}
\varphi_{s-s}(k,t)= \Psi(k e^{\beta t}),\; \beta\in \R.
\ee
Note that it corresponds to the distribution function~\eqref{2.10},
where $ a =- \beta $. The parameter $ \beta $ will be defined
below.

Then we pass to self-similar variables in~\eqref{2.15} by
substitution
\be \label{2.17}
\varphi(k,t)=\tilde{\varphi}(\tilde{k},t),
\quad
\tilde{k}=ke^{\beta t},
\ee
and obtain omitting tildes
\be
\label{2.18}
\partial_t \vphi  +\vphi +
\left(A_{\beta} k\right)\cdot \partial_k \vphi
=\Gamma\big[\vphi\big],
\quad A_{\beta}=A +{\beta}I.
\ee
It is clear that the self-similar solution~\eqref{2.16} of
Eq.~\eqref{2.15} becomes a stationary solution for
Eq.~\eqref{2.18}. The differential form of the stationary solution
is obvious from~\eqref{2.18}. Its integral form can be obtained at
 the formal level from the operator identity
\be\label{2.19}
\int\limits_{0}^{\infty} dt e^{-t \left( 1+ \hat{D}\right)
}=\left( 1+ \hat{D}\right)^{-1},
\ee
where $ \hat{D} $ is an abstract operator.
We refer to~\cite{BNV} for conditions of equivalence  of these
integral and differential forms of equation for $  \Psi(k) $. The integral equation reads~\cite{BNV}
\be \label{2.20}
 \Psi(k)=\int \limits_0^{\infty} dt E_{\beta}(t) \Gamma[\Psi(k)],
\ee
where $ \Gamma[\Psi(k)] = \mathcal{I}^{+}(\Psi,\Psi) $
is given in~\eqref{2.13},
\be \label{2.21}
 E_{\beta}(t)=\exp \left[ -t\big( 1+ A_{\beta} k\cdot \partial_k\big)\right].
 \ee
 It is easy to  see that the action of the operator  $ E_{\beta}(t)$ on any function $ \vphi(k) $ is given by formula
 \be \label{2.22}
 E_{\beta}(t)\vphi(k) =e^{-t} \vphi
 \left[e^{-{\beta}t} \big( e^{-t A}  k \big)\right].
 \ee
 The equation~\eqref{2.20} will be solved below with all
 necessary estimates. We begin in the next section with definition of $ \beta $ and some preliminary estimates.

\section{Eigenvalue problem for matrices}
\label{sec:3}
We can apply the operator
$ \big( 1+ A_{\beta} v\cdot \partial_v\big) $
to the equation~\eqref{2.20} and obtain the equation for
$ \Psi(k) $ in differential form
(see also~\eqref{2.18})
\be\label{3.1}
\big(1+{\beta} k\cdot \partial_k + A k\cdot \partial_k
\big)\Psi(k)=
\Gamma[\Psi](k).
\ee

It is always assumed below that $ \Psi(k) $ is a characteristic
 function (the Fourier transform of a probability measure in
 $ \R^{d} $) and have the following asymptotic behaviour for
 small $ \vert k \vert $:
\be\label{3.2}
\Psi(k)=1 -\frac{1}{2}B:\,k\otimes k +
O \left(|k|^{p} \right)
\ee
for some $  2 < p \leq 4 $.
The notation
$ B=\{b_{ij}; \; i,j=1, \cdots , d\} $
  is used for symmetric positively defined
matrix. We also denote for brevity
\[
B:\,k\otimes k = \sum_{i,j=1}^{n}b_{ij} k_{i}k_{j}.
\]
The formula~\eqref{3.2} means that the corresponding distribution function,
 i.e. the inverse Fourier transform of
$\Psi(k)  $, has finite moments of the order $ 2+ \varepsilon $,
$ \varepsilon > 0 $ (see~\cite{BNV} for details). It can be shown
that the matrix $ B $ and the parameter $ \beta $ satisfy the
following equation (see Eq.~\eqref{2.7} in~\cite{BNV}):
\be\label{3.3}
\beta B + \theta \left( B - \frac{\mathrm{Tr }B}{d}I\right) +
\langle B A\rangle =0,
\ee
where
\be\label{3.4}
\begin{split}
\theta=\dfrac{q d}{4(d-1)},
\quad
q=\int_{S^{d-1}} dn g(\omega \cdot n)
[1-(\omega \cdot n)^{2}],
\quad
\omega \in S^{d-1};
\\
\mathrm{Tr }B = \sum_{i=1}^{d} b_{ii},
\quad
\langle B A\rangle = \frac{1}{2}[ B A+ (B A)^{T}],
\end{split}
\ee
here the upper index $ T $ denotes the transposed matrix.
This equation can be easily obtained by substitution
of~\eqref{3.2} into Eq.~\eqref{3.1}.
We are interested in solution $( \beta, B) $ of the eigenvalue problem~\eqref{3.3} such  that the eigenvalue $ \beta  $
has the largest (as compared to other eigenvalues) real part.
In addition, the real symmetric matrix $ B  $
must have only positive eigenvalues. The existence of such solution  $( \beta, B )$ was proved in Lemma 7.3 in~\cite{BNV}
under assumption that $ \|A\| \leq \varepsilon_{0} $
for sufficiently small $ \varepsilon_{0} > 0 $,
where
\be\label{3.5}
\|A\|= \sup_{ |k| = 1 } | A k |,
\quad
k \in \R^{d}.
\ee
No estimates of $  \varepsilon_{0} $ was given in~\cite{BNV}.
Our aim in this paper is to fill this gap and to show that main
results of that paper remain valid for moderately small
values of $  \varepsilon_{0} $.
\section{Main results}
\label{sec:4}
We begin with the eigenvalue problem from Sect.~\ref{sec:3}.
\begin{lemma}
\label{lem:1}
For any real $ (d \times d) $ matrix $ A $ satisfying the
condition
\be
\label{4.1}
\|A\| < \frac{\theta}{6},
\quad
\theta=\frac{q d}{4 (d-1)},
 \ee
 in the notation of Eq.~\eqref{3.4}, the eigenvalue
  problem~\eqref{3.3} has a unique solution $( \beta, B ) $
 such that
 \begin{enumerate}
 \item{$  \beta $ has the largest real part among all
 eigenvalues of the problem~\eqref{3.3} and}
 \item{the symmetric  ($ d \times d $)-matrix $ B $ is
 normalized by condition $ \mathrm{Tr }B =d $.}
 \end{enumerate}
This solution can be represented by power series
\begin{equation*}
\beta= \theta \sum_{n=1}^{\infty}\beta_{n} \varepsilon^{n},
\quad
B =\sum_{n=1}^{\infty} B_n \varepsilon^{n},
  \quad
  \varepsilon=\frac{ \|A\|}{\theta},
  \quad B _0=I,
\end{equation*}
convergent for $ |\varepsilon| \leq 1/6 $.
The eigenvalue $ \beta $  is real and simple. The matrix
$ B $ is positive-definite. The following estimates are valid
under condition~\eqref{4.1}
\begin{equation*}
|\beta| < 2 \|A\|,
\quad
\beta - \Re \beta' \geq \theta - 5 \|A\|,
\quad
\| B - I \|  < 1,
\end{equation*}
where $ \beta' \neq \beta $ is any other eigenvalue of
the  problem~\eqref{3.3}.
\end{lemma}
This lemma defines  the parameter $ \beta $ in the
equation~\eqref{2.20} for the self-similar profile $ \Psi(k) $.
It also defines the matrix $ B $ related to behaviour
of $ \Psi(k) $ for small $ |k| $ (see Eq.~\eqref{3.2}).
The function $ \Psi(k) $ can be constructed by iterations
of the integral operator in~\eqref{2.20}. A convenient initial
approximation is given by function
\begin{equation*}
\Psi_0(k) = \exp \left( -\frac{1}{2} B:\,k\otimes k\right),
\quad
k \in  \R^{d}.
\end{equation*}
We just formulate the final result.
\begin{theorem}
\label{th:1}
Consider the integral equation~\eqref{2.20} and assume that
\be
\label{4.2}
\|A\|< \frac{q}{24},
\quad
q=\int_{S^{d-1}} dn g(\omega \cdot n) [1 - (\omega \cdot n)^{2}],
\quad
\omega \in  S^{d-1}.
\ee
It is also assumed that the solution $ \Psi(k) $ of that equation
has asymptotic behaviour for small $ |k| $ in accordance
with Eq.~\eqref{3.2} for some $  p \in (2,4] $.

(i) Then the symmetric matrix $ B  $ in~\eqref{3.2}
normalized by condition $ \mathrm{Tr }B = d $
and the parameter $\beta   $ in~\eqref{2.20} coincide with the
solution $ ( \beta, B) $  of eigenvalue problem~\eqref{3.3}
constructed in Lemma~\ref{lem:1}.

(ii) For $ \beta $ and  $ B $ from the item (i) there is a unique
characteristic function  $ \Psi(k) $ that solves
equations~\eqref{2.20}--\eqref{2.22} and satisfies the
asymptotic formula~\eqref{3.2}. The correct value of $ p $
in~\eqref{3.2} is $ p=4 $ provided that condition~\eqref{4.2}
holds.
\end{theorem}

We also consider the initial value problem for characteristic
function $ \vphi(k,t) $ in self-similar coordinates~\eqref{2.15}.
The problem reads (see Eq.~\eqref{2.18})
\be
\label{4.3}
\vphi_t + A_\beta k \cdot \vphi_k + \vphi=
\Gamma(\vphi),
\quad
\vphi|_{t=0}=\vphi_0(k),
\quad
k \in \R^{d},
\ee
where $ A_\beta = A + \beta I $, $ \beta \in \R $.
It is assumed that
\be
\label{4.4}
\left| \vphi_0(k) - \left(1 - \frac{1}{2} G_0:
k \otimes k \right)  \right| \leq C_0 |k|^{4},
\quad C_0= \mathrm{const.},
\quad
k \in \R^d,
\ee
in the notation analogous to Eq.~\eqref{3.2}.
Then it is known from~\cite{BNV} that there exists a unique
characteristic function $ \vphi(k,t) $ that solves the
problem~\eqref{4.3} and satisfies the condition
\begin{equation*}
\left| \vphi(k,t) - \left(1 - \frac{1}{2} G(t):
k \otimes k \right)  \right| \leq C_1 |k|^{4},
\quad
C_1= \mathrm{const.},
\quad
k \in \R^{d},
\end{equation*}
where $ G(t)$ is a time-dependent symmetric  $ (d \times d) $
matrix that solves the problem
\be
\label{4.5}
\frac{1}{2} G_t +  \beta G + \theta \left(G - \frac{\mathrm{Tr} G}{d} I \right)  +   \left\langle G A\right\rangle= 0,
 \quad
 G|_{t=0}= G_0,
\ee
in the notation of Eqs.~\eqref{3.4} with $ B=G $.

Hence, the matrix $ G(t) $ satisfies the linear ODE with
constant coefficients. We assume that the parameter $ \beta $
in~\eqref{4.3} is the eigenvalue from Lemma~\ref{lem:1}.
Then the matrix $ B $ from Lemma~\ref{lem:1} is a stationary
solution of Eqs.~\eqref{4.3}. We can prove that
\begin{equation*}
G(t) = c^{2 }B + O[\exp(-q t/ 12)],
\quad
t \geq 0,
\end{equation*}
for some constant $ c > 0 $ provided the condition~\eqref{4.2}
is satisfied. The constant $ c $ depends on $ G_0 $.
The following statement shows the asymptotic role of the self-similar profile $ \Psi(k) $.
\begin{theorem}
\label{th:2}
Let $ \vphi(k,t)  $ be a solution of the problem~\eqref{4.3},
where $ \|A\| < q /24 $  and $ \vphi_0(k) $ is a characteristic
function satisfying~\eqref{4.4}. Let the parameter $ \beta $
in~\eqref{4.3} and the function $   \Psi(k) $ be the same as
 in Theorem~\ref{th:1}. Then there exist two constants
 $  c>0 $ and $ C> 0 $ such that
\begin{equation*}
| \vphi(k,t) - \Psi(ck) | \leq C (|k|^{2} + |k|^{4})
e^{- \mu t},\;   \;
\mu= \frac{q - 24 \|A\| }{16},
\;
k \in \R^{d}, \; t \ge 0.
\end{equation*}
\end{theorem}
The above results can be expressed in terms of distribution
functions $ f(v,t) $ in the same way as in~\cite{BNV}.
\section{Conclusions}
We have considered the modified spatially homogeneous
Maxwell -- Boltzmann equation~\eqref{2.1}. The equation
contains an additional force term $\mathrm{div} Av f  $, where
$ v \in \R^{d} $, $ A $ is an arbitrary constant
$ (d \times d) $-matrix. Applications of this equation are
connected with well-known homoenergetic solutions to the
spatially inhomogeneous Boltzmann equation studied by many
authors since 1950s. The self-similar solutions and related questions for equation~\eqref{2.1} were recently considered
in detail in~\cite{BNV} be using  the Fourier transform
and some properties of the Boltzmann collision operator
 in the Fourier
 representation~\cite{BCG}. Main result of~\cite{BNV}
 were obtained under assumption of "sufficiently small
 norm of $ A $" in~\eqref{2.1} without explicit estimates
 of this "smallness". Our aim in this paper was to fill
 this gap and to prove that most of the results related
 to self-similar solutions remain valid for moderately
 small matrices $ A $ with norm $ \|A\|= O(10^{-1}) $
in dimensionless units. This is important for applications
because it shows boundaries for the approach based on the
 perturbation theory. The main results of the paper are
 formulated in Theorems~\ref{th:1} and~\ref{th:2} from
Sect.~\ref{sec:4}. These theorems extend the corresponding
results of~\cite{BNV} to moderate values of $ \|$ $ A\| $.
The main idea of proofs of new estimates is  based on
detailed study of the eigenvalue problem ~\eqref{3.3}, see
Lemma~\ref{lem:1} from Sect.~\ref{sec:4}. A by-product
result is the proof of existence of the bounded fourth moment of the self-similar profile for moderate values
of $ \|$ $ A\| $. The question of existence of \emph{all }
moments for the self-similar profile $ F(v) $ remains
open even in the class of arbitrarily small norm of
$ A $.
\begin{acknowledgement}
The work is dedicated to memory of Carlo Cercignani.
The research was supported by Russian Science Foundation
 Grant No. 18-11-00238-$\Pi$. The author thanks Alessia Nota and
Juan Velazquez for valuable discussions.
\end{acknowledgement}

\end{document}